\documentclass[reqno,11pt,oneside,a4paper]{amsart}

\usepackage{amsthm}
\usepackage{amscd}
\usepackage{amsmath,amssymb}
\usepackage{algorithm}
\usepackage[noend]{algorithmic}
\usepackage[format=hang]{caption}
\usepackage{url}
\usepackage{ulem}
\usepackage{csquotes}
\usepackage{multirow}
\usepackage{arydshln}
\usepackage{graphicx}
\usepackage{hyperref}
\usepackage{tikz}
\usetikzlibrary{trees}

\MakeOuterQuote{"}
\newcommand{\algorithmicbreak}{\textbf{break}}
\newcommand{\BREAK}{\STATE \algorithmicbreak}

\theoremstyle {plain}
\newtheorem {thm}{Theorem}[section]
\newtheorem {prop}[thm]{Proposition}
\theoremstyle {definition}
\newtheorem {defn}[thm]{Definition}
\newtheorem {rem}[thm]{Remark}
\newtheorem {notation}[thm]{Notation}

\DeclareMathOperator{\Ext}{Ext}
\DeclareMathOperator{\Image}{im}
\DeclareMathOperator{\Mon}{Mon}
\DeclareMathOperator{\LM}{LM}
\DeclareMathOperator{\LC}{LC}
\DeclareMathOperator{\LT}{LT}
\DeclareMathOperator{\Lead}{L}
\DeclareMathOperator{\tail}{tail}
\DeclareMathOperator{\lcm}{LCM}
\DeclareMathOperator{\Syz}{Syz}
\DeclareMathOperator{\LOT}{LOT}

\newcommand{\PP}{{\mathbb P}}
\newcommand{\Q}{{\mathbb Q}}
\newcommand{\lot}[1]{\smash{\uwave{#1}}}
\newcommand{\mcL}{{\mathcal{L}}}
\newcommand{\gen}{\ensuremath{\mathrm{e}}}
\newcommand{\Singular}{\textsc{Singular}}

\newcommand{\SyzSchreyer}{\textsc{SyzSchreyer}}
\newcommand{\LeadSyz}    {\textsc{LeadSyz}}
\newcommand{\SyzLift}    {\textsc{SyzLift}}
\newcommand{\Lift}       {\textsc{Lift}}
\newcommand{\LiftReduce} {\textsc{LiftReduce}}
\newcommand{\LiftHybrid} {\textsc{LiftHybrid}}
\newcommand{\LiftTree}   {\textsc{LiftTree}}
\newcommand{\LiftSubtree}{\textsc{LiftSubtree}}

\tikzset{
    syzdiagram/.style={
        baseline=(current bounding box.north),
        grow=down,
        level distance=20mm,
        every node/.style={rectangle, draw=black, thick, align=left},
        edge from parent fork down,
        edge from parent/.style={->, thick, draw=black}
    }
}
\tikzset{
    rightof/.style args={#1}{
        at=(#1.east),
        anchor=west,
        shift=(0:3mm)
    }
}

\begin{document}

\bibliographystyle{alpha}

\title{Refined Algorithms to Compute Syzygies}

\author{Bur\c{c}in Er\"ocal}
\address{Bur\c{c}in Er\"ocal\\
Department of Mathematics\\
University of Kaiserslautern\\
Erwin-Schr\"odinger-Str.\\
67663 Kaiserslautern\\
Germany}
\email{burcin@mathematik.uni-kl.de}

\author{Oleksandr Motsak}
\address{Oleksandr Motsak\\
Department of Mathematics\\
University of Kaiserslautern\\
Erwin-Schr\"odinger-Str.\\
67663 Kaiserslautern\\
Germany}
\email{motsak@mathematik.uni-kl.de}
\urladdr{http://www.mathematik.uni-kl.de/$\sim$motsak}

\author{Frank-Olaf Schreyer}
\address{Frank-Olaf Schreyer\\
Mathematik und Informatik\\
Universit\"at des Saarlandes\\
Campus E2 4\\
66123 Saarbr\"ucken\\
Germany}
\email{schreyer@math.uni-sb.de}
\urladdr{http://www.math.uni-sb.de/ag/schreyer}

\author{Andreas Steenpa\ss}
\address{Andreas Steenpa\ss\\
Department of Mathematics\\
University of Kaiserslautern\\
Erwin-Schr\"odinger-Str.\\
67663 Kaiserslautern\\
Germany}
\email{steenpass@mathematik.uni-kl.de}


\date{\today}

\begin{abstract}
Based on Schreyer's algorithm \cite{Schreyer1, Schreyer2, BS}, we present two
refined algorithms for the computation of syzygies. The two main ideas of the
first algorithm, called \hyperref[alg:LiftHybrid]{\LiftHybrid{}}, are the
following: First, we may leave out certain terms of module elements during the
computation which do not contribute to the result. These terms are called
"lower order terms", see Definition~\ref{def:lower_order_term}. Second, we do
not need to order the remaining terms of these module elements during the
computation. This significantly reduces the number of monomial comparisons for
the arithmetic operations. For the second algorithm, called
\hyperref[alg:LiftTree]{\LiftTree{}}, we additionally cache some partial
results and reuse them at the remaining steps.
\end{abstract}

\maketitle

\section{Introduction}

Computing syzygies, that is, a free resolution
\[
0 \longrightarrow F_n \longrightarrow \ldots \longrightarrow F_1
\longrightarrow F_0 \longrightarrow M \longrightarrow 0
\]
of a module $M$ over a polynomial ring $R = K[x_1, \ldots, x_n]$, is one of the
fundamental tasks in constructive module theory, needed for example for the
computation of the modules $\Ext^i_S(M, N)$ for a further $R$-module $N$. An
algorithm for computing the resolution, starting from a Gr\"obner basis of the
image of the presentation matrix $\Image(F_1 \to F_0)$, was given by Schreyer
\cite{Schreyer1}, see also \cite[Theorem 15.10, Corollary 15.11]{Eis95} or
\cite[Corollary 1.11]{BS}.

If $M$ is a graded module over the polynomial ring $R$ with its standard
grading, then there exists a minimal free resolution which is uniquely
determined up to isomorphism. While the computation of the Gr\"obner basis of
the presentation matrix is still feasible in many examples, the computation of
the minimal free resolution might be out of reach. However, the computation of
a non-minimal resolution is typically much cheaper, and for many applications,
such as the computation of a single $\Ext$ module, good enough.

In this paper we describe a refined version of Schreyer's algorithm, which
utilizes the full strength of Schreyer's Theorem \cite[Theorem 15.10]{Eis95}.
The basic idea, which can be already found in \cite{Schreyer2}, is to ignore
lower order terms in the computation of the generators of the next syzygy
module. This is possible since these terms will cancel each other anyway.

Our implementation gives a considerable speed-up, in many cases even if we
additionally minimize the resolution. The computation of the minimal Betti
numbers from a non-minimal resolution is typically much faster than minimizing
the whole resolution. And, for large examples, the computation of a non-minimal
resolution using our method is again faster than deriving the minimal Betti
numbers from it. Note that in many cases, however, even a non-minimal
resolution suffices to deduce geometric information, see, for example,
Remark~\ref{rem:pcnc_g16}.

Moreover, our findings suggest that starting, say, from a Gr\"obner basis over
$\Q$, the computation of a non-minimal resolution using floating point numbers
as coefficients might be numerically stable. However, this is a topic of future
work.

The paper is organized along the following lines. In
Section~\ref{sec:preliminary}, we introduce some basic terminology. The induced
monomial ordering, Schreyer's Theorem, and the corresponding algorithm are
discussed in Section~\ref{sec:schreyer_alg}. Based on an analysis of this
algorithm, we present the two new algorithms in
Section~\ref{sec:new_algorithms}. A detailed example is given in
Section~\ref{sec:example}. In Section~\ref{sec:timings}, we illustrate our
implementation in a number of examples. One series of examples consists of
Artinian graded Gorenstein algebras which we regard as an appropriate family of
examples to test any syzygy algorithm since we can vary the number of
variables, the degree, and the sparseness. Further series are nodal canonical
and nodal Prym canonical curves. For the examples, we always work over a finite
ground field.

\section{Preliminaries}\label{sec:preliminary}

Throughout this article, let $K$ be a field, and let $R := K[x_1, \ldots, x_n]$
be the polynomial ring in $n$ variables over $K$. We denote the monoid of
monomials in $x_1, \ldots, x_n$ by $\Mon(x_1, \ldots, x_n)$.

We briefly recall some terminology for dealing with $R$-module syzygies and
their computation.

\begin{defn}
Let $F := R^r$ be the free $R$-module of rank $r$, and let
$\gen_1, \ldots, \gen_r$ be the canonical basis of $F$.
\begin{enumerate}
\item
A \emph{monomial} in $F$ is the product of an element in
$\Mon(x_1, \ldots, x_n)$ with a basis element $\gen_i$. The set of monomials
in $F$ is denoted by $\Mon(F)$.

\item Accordingly, a \emph{term} in $F$ is the product of a monomial in $F$
with a scalar in $K$.

\item
A monomial $m_1 \gen_i$ \emph{divides} a monomial $m_2 \gen_j$ if $i = j$ and
$m_1$ divides $m_2$; in this case, the \emph{quotient}
$m_2 \gen_j / m_1 \gen_i$ is defined as $m_2 / m_1 \in \Mon(x_1, \ldots, x_n)$.
We also say that $m_1 \in \Mon(x_1, \ldots, x_n)$ divides $m_2 \gen_j$ if $m_1$
divides $m_2$, and in this case $m_2 \gen_j / m_1$ is defined to be
$(m_2 / m_1) \gen_j \in \Mon(F)$.

\item The \emph{least common multiple} of two monomials
$m_1 \gen_i, m_2 \gen_j \in F$ is
\[
\lcm(m_1 \gen_i, m_2 \gen_j)
:= \left\{ \begin{array}{ll}
             \lcm(m_1, m_2) \gen_i, &\text{ if } i = j, \\
             0,                            &\text{ otherwise.}
           \end{array}
   \right.
\]

\item
A \emph{monomial ordering} on $F$ is a total ordering $\succ$ on $\Mon(F)$ such
that if $m_1 \gen_i$ and $m_2 \gen_j$ are monomials in $F$, and $m$ is a
monomial in $R$, then
\[
m_1 \gen_i \succ m_2 \gen_j
\;\Longrightarrow\; (m \cdot m_1) \gen_i \succ (m \cdot m_2) \gen_j \,.
\]
In this article, we require in addition that
\[
m_1 \gen_i \succ m_2 \gen_i
\;\Longleftrightarrow\; m_1 \gen_j \succ m_2 \gen_j
\text{ for all } i, j \,.
\]

\item
Let $\succ$ be a monomial ordering on $F$, let $f \in F \setminus \{0\}$ be an
element of $F$, and let $f = cm \gen_i + f^*$ be the unique decomposition of
$f$ with $c \in K \setminus \{0\}$, $m \gen_i \in \Mon(F)$, and
$m \gen_i > m^* \gen_j$ for any non-zero term $c^* m^* \gen_j$ of $f^*$.
We define the \emph{leading monomial}, the \emph{leading coefficient}, the
\emph{leading term}, and the \emph{tail} of $f$ as
\begin{align*}
\LM(f)   &:= m \gen_i  \,, \\
\LC(f)   &:= c         \,, \\
\LT(f)   &:= cm \gen_i \,, \\
\tail(f) &:= f-\LT(f)  \,,
\end{align*}
respectively.

\item For any subset $S \subset F$, we call
\[
\Lead(S) := \langle \LM(f) \mid f \in S \setminus \{0\} \rangle_R \subset F
\]
the \emph{leading module} of $S$.
\end{enumerate}
\end{defn}

\begin{rem}
Let $\succ$ be a monomial ordering on the free $R$-module $F := R^r$ as defined
above. Then there is a unique monomial ordering $>$ on $R$ which is compatible
with $\succ$ in the obvious way, and we say that $\succ$ is \emph{global} if
$>$ is global. In this article, all monomial orderings are supposed to be
global.
\end{rem}

\begin{defn}\label{def:syzygy}
Let $N$ be an $R$-module, let $G := \{f_1, \ldots, f_r\} \subset N$ be a finite
subset of $N$, and let $F := R^r$ be the free $R$-module of rank $r$ as above.
Consider the homomorphism
\begin{align*}
\psi_G: F &\rightarrow N \,, \\
\gen_i &\mapsto f_i \,.
\end{align*}
A \emph{syzygy} of $G = \{f_1, \ldots, f_r\}$ is an element of $\ker{\psi_G}$.
We call $\ker{\psi_G}$ the (first) \emph{syzygy module} of $G$, written
\[
\Syz(G) := \ker{\psi_G}.
\]
\end{defn}

\begin{defn}\label{def:free_resolution}
Let $M$ be an $R$-module. A \emph{free resolution} of $M$ is an exact sequence
\[
\mathcal{F}: \quad
\ldots \longrightarrow F_{i+1}
\overset{\phi_{i+1}}{\longrightarrow} F_i
\overset{\phi_i}{\longrightarrow} F_{i-1}
\longrightarrow \ldots
\longrightarrow F_1
\overset{\phi_1}{\longrightarrow} F_0
\longrightarrow M
\longrightarrow 0
\]
with free $R$-modules $F_i$, $i \in \mathbb{N}$.
\end{defn}

\begin{rem}\label{rem:free_resolution}
Let the notation be as in Definitions~\ref{def:syzygy} and
\ref{def:free_resolution}. In this article, we only consider the case where $N$
is a free module over the polynomial ring $R$ and where we wish to construct a
free resolution of $M = N / \langle G \rangle_R$. For this, with notation as in
Definition~\ref{def:syzygy}, set $F_0 := N$, $F_1 := F$, and
$\phi_1 := \psi_G$. Now, starting with $G_1 := G$, let $G_{i+1}$ be a finite
set of generators for $\Syz(G_i)$ and, inductively, define $\phi_i$ to be the
map $\psi_{G_i}$ for $i \in \mathbb{N} \setminus \{0\}$. We then have
$\Syz(G_i) = \ker \phi_i$, that is, $\mathcal{F}$ is obtained by repeatingly
computing the syzygies of finite subsets of free $R$-modules.
\end{rem}

\begin{defn}\label{def:std_representation}
Let $F_0 := R^s$ be the free $R$-module of rank $s$, let $>$ be a monomial
ordering on $F_0$, and let
$G := \{f_1, \ldots, f_r\} \subset F_0 \setminus \{0\}$ be a set of non-zero
vectors in $F_0$.

\begin{enumerate}
\item
We define $m_{ji}$ as
\[
m_{ji} := \frac{\lcm(\LM(f_j), \LM(f_i))}{\LT(f_i)} \in R \,.
\]

\item\label{item:m_ij}
For $i,j \in \{1, \ldots, r\}$, we define the \emph{S-vector} of $f_i$ and
$f_j$ as
\[
S(f_i, f_j) := m_{ji} f_i - m_{ij} f_j \in \langle G \rangle_R \subset F_0 \,.
\]

\item\label{item:std_representation}
For $g \in F_0$, we call an expression
\[
g = g_1 f_1 + \ldots + g_r f_r + h
\]
with $g_i \in R$ and $h \in F_0$ a \emph{standard representation} for $g$ with
remainder $h$ (and w.r.t.\@ $G$ and $>$) if the following conditions are
satisfied:
\begin{enumerate}
\item
$\LM(g) \geq \LM(g_i f_i)$ for all $i = 1, \ldots, r$ whenever both $g$ and
$g_i f_i$ are non-zero.

\item
If $h$ is non-zero, then $\LT(h)$ is not divisible by any $\LT(f_i)$.
\end{enumerate}
\end{enumerate}
\end{defn}

\begin{rem}
Standard representations can be computed by multivariate division with
remainder. With notation as above, let now $G$ be a Gr\"obner basis, and let
$g$ be an element of $\langle G \rangle_R$. In this case, the remainder $h$ is
zero by Buchberger's criterion for Gr\"obner bases. For S-vectors of elements
of $G$, each standard representation
\[
S(f_i, f_j) = m_{ji} f_i - m_{ij} f_j
= g_1^{(ij)} f_1 + \ldots + g_r^{(ij)} f_r
\]
yields an element
$m_{ji} \gen_i - m_{ij} \gen_j
- \big( g_1^{(ij)} \gen_1 + \ldots + g_r^{(ij)} \gen_r \big) \in \Syz(G)$.
\end{rem}

This gives one possibility to compute syzygies which we will now discuss in
detail.

\section{Schreyer's Syzygy Algorithm}\label{sec:schreyer_alg}

\subsection{The Induced Ordering}

\begin{defn}
Given a monomial ordering $>$ on $F_0 := R^s$ and a set of non-zero vectors
$G := \{f_1, \ldots, f_r\} \subset F_0 \setminus \{0\}$, we define the
\emph{induced ordering} on $F_1 := R^r$ (w.r.t.\@ $>$ and $G$) as the monomial
ordering $\succ$ given by
\begin{align*}
m_1 \gen_i \succ m_2 \gen_j \; :\Leftrightarrow \;
&\LT(m_1 f_i) > \LT(m_2 f_j) \\
&\text{ or } (\LT(m_1 f_i) = \LT(m_2 f_j) \text{ and } i > j)
\end{align*}
for all monomials $m_1, m_2 \in \Mon(x_1 \ldots, x_n)$, and for all basis
elements $\gen_i, \gen_j \in F_1$.
\end{defn}

This definition implies that both $>$ and $\succ$ yield the same ordering on
$R$ if restricted to one component.

Monomial comparisons w.r.t.\@ induced orderings are computationally expensive
and should therefore be avoided in practice. This holds in particular in the
case of chains $(\succ_i)_{i=1,\ldots,k}$ of orderings with $\succ_{i+1}$
induced by $\succ_i$ which appear in the computation of free resolutions.

\subsection{Schreyer's Theorem}

\begin{thm}[{\cite[Corollary~1.11]{BS}}]\label{thm:Schreyer}\ \\
Let $G = \{f_1, \ldots, f_r\} \subset F_0 := R^s$ be a Gr\"obner basis w.r.t.\@
a monomial ordering~$>$ on $F_0$. For each pair $(f_i, f_j)$ with
$i, j \in \{1, \ldots, r\}$, let
\[
S(f_i, f_j) = m_{ji} f_i - m_{ij} f_j
= g_1^{(ij)} f_1 + \ldots + g_r^{(ij)} f_r
\]
be a standard representation of the corresponding S-vector.
Then the relations
\[
m_{ji} \gen_i - m_{ij} \gen_j
- \left( g_1^{(ij)} \gen_1 + \ldots + g_r^{(ij)} \gen_r \right)
\in F_1 := R^r
\]
form a Gr\"obner basis of $\Syz(G)$ w.r.t.\@ the monomial ordering on $F_1$
induced by $>$ and $G$. In particular, these relations generate the syzygy
module $\Syz(G)$.
\end{thm}

Based on this theorem, there is an obvious algorithm for the computation of
syzygy modules: Given a Gr\"obner basis $G$ as above, it suffices to compute
standard representations for all S-vectors $S(f_i, f_j)$ by division with
remainder.

Of course, one can do much better. Since $S(f_i, f_j) = -S(f_j, f_i)$, it is
sufficient to consider those pairs $(f_i, f_j)$ with $j < i$. It is well-known
that even more pairs can be left out using the following notation
(cf.\@ \cite{BS}):

\begin{notation}\label{notation:M_i}
Let $F_0 := R^s$ be the free $R$-module of rank $s$, and let
$G := \{f_1, \ldots, f_r\} \subset F_0 \setminus \{0\}$ be a set of non-zero
vectors in $F_0$. For $i = 2, \ldots, r$, we define the monomial ideal $M_i$ as
\[
M_i := \langle \LT(f_1), \ldots, \LT(f_{i-1}) \rangle
    : \langle \LT(f_i) \rangle \subseteq R \,.
\]
\end{notation}

\begin{rem}
Recall that if $N_1$ and $N_2$ are submodules of an $R$-module $M$, then the
module quotient $N_1 : N_2$ is defined to be the ideal
\[
N_1 : N_2 := \{ a \in R \mid an \in N_1 \text{ for all } n \in N_2 \}
\subseteq R \,.
\]
In particular, in the situation of Notation~\ref{notation:M_i}, we have
$\langle m_1 \gen'_i \rangle : \langle m_2 \gen'_j \rangle = 0$ for any two
monomials $m_1, m_2 \in R$ and any two basis elements $\gen'_i$, $\gen'_j$ of
$F_0$ with $i \neq j$.
\end{rem}

\begin{prop}[{\cite[Theorem~1.5]{BS}}]
Let $G = \{f_1, \ldots, f_r\} \subset F_0$ be as in Theorem~\ref{thm:Schreyer}.
For each $i = 2, \ldots, r$, and for each minimal generator $x^\alpha$ of the
monomial ideal $M_i \subset R$, let $j = j(i, \alpha) < i$ be an index such
that $m_{ji}$ divides $x^\alpha$. Then it is sufficient in
Theorem~\ref{thm:Schreyer} to consider only the corresponding pairs
$(f_i, f_j)$.
\end{prop}

Taking this proposition into account, we get Algorithm~\ref{alg:SyzSchreyer}
below.

\begin{algorithm}[htb]
\caption{\SyzSchreyer}\label{alg:SyzSchreyer}
\begin{algorithmic}[1]
\REQUIRE A Gr\"obner basis $G = \{f_1, \ldots, f_r\} \subset F_0 := R^s$
w.r.t.\@ some monomial ordering $>$
\ENSURE A Gr\"obner basis of $\Syz(G) \subset F_1 := R^r$ w.r.t.\@ the monomial
ordering induced by $>$ and $G$
\vspace{0.2cm}

\STATE $S := \varnothing$
\FOR{$i = 2, \ldots, r$}
  \FOR{each minimal generator $x^\alpha$ of the monomial ideal $M_i$}
    \STATE choose an index $j < i$ such that $m_{ji}$ divides $x^\alpha$
    \STATE $h := S(f_i, f_j) = m_{ji} f_i - m_{ij} f_j \in F_0$
    \STATE $s := m_{ji} \gen_i - m_{ij} \gen_j \in F_1$
    \WHILE{$h \neq 0$}
      \STATE choose an index $\lambda$ such that $\LT(f_\lambda)$ divides
          $\LT(h)$
      \STATE $h := h - \frac{\LT(h)}{\LT(f_\lambda)} f_\lambda$
      \STATE $s := s - \frac{\LT(h)}{\LT(f_\lambda)} \gen_\lambda$
    \ENDWHILE
    \STATE $S := S \cup \{ s \}$
  \ENDFOR
\ENDFOR
\RETURN $S$
\end{algorithmic}
\end{algorithm}

\subsection{Schreyer Frame}\label{ssec:schreyer_frame}

The leading module of the syzygy module will serve as a starting point for the
algorithms which we propose in Section~\ref{sec:new_algorithms}. Its
computation is based on the following observation.

\begin{rem}\label{rem:lead_syz}
With notation as in Theorem~\ref{thm:Schreyer},
$g_1^{(ij)} f_1 + \ldots + g_r^{(ij)} f_r$ is a standard representation of the
S-vector $S(f_i, f_j)$, and therefore we have
$\LM(m_{ji} f_i) = \LM(m_{ij} f_j) > \LM\Bigl(g^{(ij)}_k f_k\Bigr)$ for all
$k = 1, \ldots, r$ with $g^{(ij)}_k \neq 0$,
cf.\@ Definition~\ref{def:std_representation}(\ref{item:std_representation}).
For $i > j$, this implies
\[
m_{ji} \gen_i \succ m_{ij} \gen_j \succ \LM\Bigl(g^{(ij)}_k\Bigr) \gen_k \,,
\]
where $\succ$ is the monomial ordering on $F_1 = R^r$ induced by $>$ and $G$.
Therefore the leading syzygy module of $G$ w.r.t.\@ $\succ$ is
\[
\Lead_\succ(\Syz(G)) = \bigoplus_{i=2, \ldots, r} M_i \, \gen_i \,.
\]
\end{rem}

Thus, for a given Gr\"obner basis $G$, the leading module of $\Syz(G)$ w.r.t.\@
the induced ordering can be easily computed by throwing away superfluous
elements, see Algorithm~\ref{alg:LeadSyz}.

\begin{algorithm}[htb]
\caption{\LeadSyz}\label{alg:LeadSyz}
\begin{algorithmic}[1]
\REQUIRE A Gr\"obner basis $G = \{f_1, \ldots, f_r\} \subset F_0 := R^s$
w.r.t.\@ some monomial ordering $>$ on $F_0$
\ENSURE A minimal set $\mcL$ of generators for the leading syzygy module
$\Lead_\succ(\Syz(G))$ of $G$ w.r.t.\@ the monomial ordering $\succ$ on
$F_1 := R^r$ induced by $>$ and $G$
\vspace{0.2cm}

\STATE $\mcL := \varnothing$
\FOR{$1 \leq j < i \leq r$}
  \STATE $t := m_{ji} \gen_i \in F_1$
  \FOR{$s \in \mcL$}
    \IF{$s \mid t$}
      \STATE $t:=0$
      \BREAK
    \ELSIF{$t \mid s$}
      \STATE $\mcL := \mcL \setminus \{s\}$
    \ENDIF
  \ENDFOR
  \IF{$t\neq0$}
    \STATE $\mcL := \mcL \cup \{t\}$
  \ENDIF
\ENDFOR
\RETURN $\mcL$
\end{algorithmic}
\end{algorithm}

In Algorithm~\ref{alg:LeadSyz}, only the leading terms of the Gr\"obner
basis $G$ contribute to the computation of the set $S$ (via the term
$m_{ji} \in R$, cf.\@
Definition~\ref{def:std_representation}(\ref{item:m_ij})). For a free
resolution as constructed in Remark~\ref{rem:free_resolution}, we can thus,
starting with the leading terms of $G$, inductively compute sets of generators
for all leading syzygy modules. The sequence of these sets of leading syzygy
terms is called a \emph{Schreyer frame} by La Scala and Stillman in \cite{LS}.

It is worth noting that the algorithm to compute a minimal free resolution by
La~Scala and Stillman is compatible with our algorithms for the computation of
syzygies in the sense that both approaches are based on the Schreyer frame and
can thus be combined.

\begin{rem}
In the computation of a free resolution, reordering the syzygies after each
step may yield smaller generators for higher syzygy modules. With notation as
in Remark~\ref{rem:free_resolution}, we expect that reordering $G_i$ w.r.t.\@
the negative degree reverse lexicographical ordering on $F_{i-1}$ before
computing $G_{i+1}$ is generally the best choice.
\end{rem}

\section{New Algorithms}\label{sec:new_algorithms}

Throughout this section, let $G := \{f_1, \ldots, f_r\} \subset F_0 := R^s$ be
a Gr\"obner basis w.r.t.\@ some monomial ordering $>$ and let $\succ$ be the
monomial ordering on $F_1 := R^r$ induced by $>$ and $G$. Furthermore, let
$\mcL$ be the minimal generating set of the monomial submodule
$\Lead_\succ(\Syz(G)) \subset F_1$. We simply write $\psi$ for the map
$\psi_G: F_1 \rightarrow F_0$ defined by $\psi_G(\gen_i) := f_i$ as in
Definition~\ref{def:syzygy}.

By Remark~\ref{rem:lead_syz}, there is a one-to-one correspondence between the
minimal generators of the monomial ideals $M_i$ and the elements of $\mcL$.
Instead of processing S-pairs, we can therefore directly start with the minimal
generating set of leading syzygy terms. This is equivalent to applying the
chain criterion for syzygies to the set of all S-pairs,
cf.~\cite[Lemma~2.5.10]{GP}.

The algorithmic idea is that each leading syzygy term $s \in \mcL$ gives rise
to a pair of indices $(i, j)$ with $s = m_{ji} \gen_i$, which, through a
standard representation of the corresponding S-vector $S(f_i, f_j)$, gives rise
to a syzygy~$\bar{s}$ of $G$ with $\LT_\succ(\bar{s}) = s$. Note that both
the pair of indices and the standard representation obtained thereof are in
general not unique.

This motivates the following definition.

\begin{defn}\label{def:lifting}
Let $s \in \Lead_\succ(\Syz(G)) \subset F_1$ be a leading syzygy term. We call
$\bar{s} \in F_1$ a \emph{lifting} of $s$ w.r.t.\@ $G$ and $\succ$ if the
following conditions hold:
\begin{enumerate}
\item
$\LT_\succ(\bar{s}) = s$, and

\item
$\bar{s} \in \Syz(G)$.
\end{enumerate}
\end{defn}

If we know how to compute such a lifting, then we can use
Algorithm~\ref{alg:SyzLift} to obtain a generating set $S$ of the syzygy
module. Since $\Lead_\succ(S)$ is equal to $\Lead_\succ(\Syz(G))$, this set is
even a Gr\"obner basis of $\Syz(G)$ w.r.t.\@ $\succ$. From the computational
point of view, Algorithm~\ref{alg:SyzSchreyer} can be regarded as the special
case of Algorithm~\ref{alg:SyzLift} where the liftings are computed by the
usual reduction. This can be reformulated as in Algorithm~\ref{alg:LiftReduce}.

\begin{algorithm}[htb]
\caption{\SyzLift}\label{alg:SyzLift}
\begin{algorithmic}[1]
\REQUIRE A Gr\"obner basis $G \subset F_0$ w.r.t.\@ $>$ and an algorithm
\Lift{} to compute, for a leading syzygy term $s \in \Lead_\succ(\Syz(G))$, a
lifting w.r.t.\@ $G$ and $\succ$ \\
(\Lift{} can be, for example, any of the three algorithms \LiftReduce{},
\LiftHybrid{}, or \LiftTree{} below.)
\ENSURE A Gr\"obner basis of $\Syz(G) \subset F_1$ w.r.t.\@ $\succ$
\vspace{0.2cm}

\STATE $\mcL := \LeadSyz(G)$
\STATE $S := \varnothing$
\FOR{$s \in \mcL$}
  \STATE $\bar{s} := \Lift(s)$
  \STATE $S := S \cup \{ \bar{s} \}$
\ENDFOR
\RETURN $S$
\end{algorithmic}
\end{algorithm}

Let us now discuss algorithms for lifting leading syzygy terms in detail.
\LiftReduce{} (Algorithm~\ref{alg:LiftReduce}) computes a lifting of a given
leading syzygy term $s \in \Lead_\succ(\Syz(G))$ via multivariate division of
the polynomial $g := \psi(s) \in \langle G \rangle \subset F_0$ by the elements
of $G$. This is computationally the same as the division of $h$ w.r.t.\@ $G$ in
the while-loop of \SyzSchreyer{} (Algorithm~\ref{alg:SyzSchreyer}). At each
step, the leading term of $g$ is reduced, and this process finally reaches
$g = 0$ since $G$ is a Gr\"obner basis.

\begin{algorithm}[htb]
\caption{\LiftReduce}\label{alg:LiftReduce}
\begin{algorithmic}[1]
\REQUIRE A Gr\"obner basis $G = \{f_1, \ldots, f_r\} \subset F_0$ w.r.t.\@ $>$
and a leading syzygy term $s \in \Lead_\succ(\Syz(G)) \subset F_1$
\ENSURE A lifting $\bar{s} \in \Syz(G) \subset F_1$ of $s$ w.r.t.\@ $G$ and
$\succ$
\vspace{0.2cm}

\STATE $g := \psi(s)$
\STATE $\bar{s} := s$
\WHILE{$g \neq 0$}
  \STATE $t := \LT(g)$
  \STATE\label{line:LiftReduce_choose}%
    choose a term $m \gen_i \in F_1$ with $m \LT(f_i) = t$ and
    $s \succ m \gen_i$
  \STATE\label{line:LiftReduce_reduction} $g := g - m f_i$
  \STATE $\bar{s} := \bar{s} - m \gen_i$
\ENDWHILE
\RETURN $\bar{s}$
\end{algorithmic}
\end{algorithm}

Let $g_1, \ldots, g_k \in F_0$ be the sequence of values which $g$ takes when
the algorithm \LiftReduce{} is applied to a leading syzygy term
$s \in \Lead_\succ(\Syz(G))$. Since we have $g_k = 0$, every single term
occurring in this sequence is eventually cancelled at one of the reduction
steps in line~\ref{line:LiftReduce_reduction}, but only the processing of the
leading terms $\LT(g_1), \ldots, \LT(g_k)$ contributes to the syzygy
$\bar{s} \in \Syz(G)$. In particular, those terms which are not divisible by
one of the leading monomials $\LM(f_i)$, $i = 1, \ldots, r$, do not contribute
to $\bar{s}$ and can therefore be left out. We use the following terminology to
refer to these terms.

\begin{defn}\label{def:lower_order_term}
Let $S \subset F_0$ be a set of vectors and let $t \in F_0$ be a term. Then $t$
is called a \emph{lower order term} w.r.t.\@ $S$ if
\[
\LM(f) \nmid t \; \text{ for all } f \in S \setminus \{0\} \,.
\]
For an element $g \in F_0$, we define $\LOT(g|S)$ to be the sum of those terms
occuring in $g$ which are of lower order w.r.t.\@ $S$.
\end{defn}

Furthermore, instead of reducing the leading term of $g$ at a given step, we
may choose any term of $g$ which is not of lower order. Taking the above
observations into account, we get the algorithm \LiftHybrid{}
(Algorithm~\ref{alg:LiftHybrid}). Note that the lower order terms which are
left out at the intermediate steps sum up to zero.

\begin{algorithm}[htb]
\caption{\LiftHybrid}\label{alg:LiftHybrid}
\begin{algorithmic}[1]
\REQUIRE A Gr\"obner basis $G = \{f_1, \ldots, f_r\} \subset F_0$ w.r.t.\@ $>$
and a leading syzygy term $s \in \Lead_\succ(\Syz(G)) \subset F_1$
\ENSURE A lifting $\bar{s} \in \Syz(G) \subset F_1$ of $s$ w.r.t.\@ $G$ and
$\succ$
\vspace{0.2cm}

\STATE $g := \psi(s)-\LOT(\psi(s)|G)$
\STATE $\bar{s} := s$
\WHILE{$g \neq 0$}
  \STATE choose a term $t$ of $g$
  \STATE choose a term $m \gen_i \in F_1$ with $m \LT(f_i) = t$ and
    $s \succ m \gen_i$
  \STATE\label{line:LiftHybrid_reduction}
    $g := g - ( m f_i - \LOT(m f_i | G) )$
  \STATE $\bar{s} := \bar{s} - m \gen_i$
\ENDWHILE
\RETURN $\bar{s}$
\end{algorithmic}
\end{algorithm}

We can even go further and consider the set $T$ of terms in $g$ rather than the
polynomial $g$ itself. In other words, we do not need to sort the terms in $g$
and we do not need to carry out the cancellations of terms which may occur in
line~\ref{line:LiftHybrid_reduction} of \LiftHybrid{}. Then each term in $T$
can be reduced independently as in \LiftTree{} (Algorithm~\ref{alg:LiftTree}).
This yields a tree structure by the recursive calls of \LiftSubtree{}
(Algorithm~\ref{alg:LiftSubtree}) for each term in $T$.

The algorithm applied at the root node of this tree, \LiftTree{}, slightly
differs from the algorithm applied at the other nodes, \LiftSubtree{}. In
\LiftTree{}, the leading term of $\psi(s)$ is included in $T$, whereas at the
other nodes, this term has been cancelled by the reduction in the previous step
and is therefore left out in \LiftSubtree{}. Because of this difference, we
need the following definition to give a proper description of the output of
\LiftSubtree{}.

\begin{defn}\label{def:subtree_lifting}
Let $s \in F_1$ be a term. We call $\hat{s} \in F_1$ a \emph{subtree lifting}
of $s$ w.r.t.\@ $G$ and $\succ$ if the following conditions hold:
\begin{enumerate}
\item
$\LT_\succ(\hat{s}) = s$, and

\item
all terms in $\tail(\psi(\hat{s})) \in F_0$ are lower order terms w.r.t.\@ $G$
and $\succ$.
\end{enumerate}
\end{defn}

\begin{algorithm}[tb]
\caption{\LiftTree}\label{alg:LiftTree}
\begin{algorithmic}[1]
\REQUIRE A Gr\"obner basis $G = \{f_1, \ldots, f_r\} \subset F_0$ w.r.t.\@ $>$
and a leading syzygy term $s \in \Lead_\succ(\Syz(G)) \subset F_1$
\ENSURE A lifting $\bar{s} \in \Syz(G) \subset F_1$ of $s$ w.r.t.\@ $G$ and
$\succ$
\vspace{0.2cm}

\STATE $g := \psi(s)$
\STATE $T :=$ set of terms in $(g-\LOT(g|G))$
\STATE $\bar{s} := s$
\FOR{all $t \in T$}
  \STATE\label{line:LiftTree_choose} choose a term $m \gen_i \in F_1$ with
    $m \LT(f_i) = t$ and $s \succ m \gen_i$
  \STATE $\bar{s} := \bar{s} - \LiftSubtree(m \gen_i)$
\ENDFOR
\RETURN $\bar{s}$
\end{algorithmic}
\end{algorithm}

\addtocounter{algorithm}{-1}
\renewcommand{\thealgorithm}{\arabic{algorithm}a}

\begin{algorithm}[tb]
\caption{\LiftSubtree}\label{alg:LiftSubtree}
\begin{algorithmic}[1]
\REQUIRE A Gr\"obner basis $G = \{f_1, \ldots, f_r\} \subset F_0$ w.r.t.\@ $>$
and a term $s \in F_1$
\ENSURE A subtree lifting $\hat{s} \in F_1$ of $s$ w.r.t.\@ $G$ and $\succ$
\vspace{0.2cm}

\STATE $g := \psi(s)-\LT(\psi(s))$
\STATE $T :=$ set of terms in $(g-\LOT(g|G))$
\STATE $\hat{s} := s$
\FOR{all $t \in T$}
  \STATE choose a term $m \gen_i \in F_1$ with $m \LT(f_i) = t$
  \STATE $\hat{s} := \hat{s} - \LiftSubtree(m \gen_i)$
\ENDFOR
\RETURN $\hat{s}$
\end{algorithmic}
\end{algorithm}

\renewcommand{\thealgorithm}{\arabic{algorithm}}

\LiftTree{} terminates when $T = \varnothing$ is reached in every branch of the
tree. One can easily check that this algorithm returns indeed a lifting of the
input by comparing it to \LiftHybrid{}.

\begin{rem}
Let $s \in \Lead_\succ(\Syz(G)) \subset F_1$ be a leading syzygy term. If
$\bar{s}$ is a lifting of $s$ w.r.t.\@ $G$ and $\succ$, then $\bar{s}$ is a
subtree lifting of $s$, but the converse statement is not true in general.

For any term $s' \in F_1$, a proper lifting of $s'$ (w.r.t.\@ $G$ and $\succ$)
exists if and only if $s'$ is leading syzygy term, that is, an element of
$\Lead_\succ(\Syz(G))$. Hence we cannot expect to find proper liftings of the
terms $m \gen_i \in F_1$ to which \LiftSubtree{} is applied in
Algorithms~\ref{alg:LiftTree} and \ref{alg:LiftSubtree}.
\end{rem}

\begin{rem}\label{rem:LiftTree_condition}
The condition $s \succ m \gen_i$ in line~\ref{line:LiftTree_choose} of
Algorithm~\ref{alg:LiftTree} is always satisfied at the analogous step in
Algorithm~\ref{alg:LiftSubtree} and thus does not need to be checked there.
\end{rem}

\begin{rem}
Since the only differences between the algorithms \LiftTree{} and
\LiftSubtree{} are the assignment of $g$ in line~1 and, as explained in
Remark~\ref{rem:LiftTree_condition} above, the condition $s \succ m \gen_i$ in
line~5, we could have merged them both into one algorithm. This can be done,
for example, by using a boolean variable which is set to true in the case of
\LiftTree{}, corresponding to the root node of the resulting tree, and which is
set to false for \LiftSubtree{}, representing the inner nodes and the leaves.
However, we think that separating the two algorithms may help to understand the
mathematical properties of the output of \LiftSubtree{} as described in
Definition~\ref{def:subtree_lifting} in contrast to the output of \LiftTree,
see Definition~\ref{def:lifting}.
\end{rem}

\begin{rem}
Let $m_1 \gen_{i_1}, \ldots, m_k \gen_{i_k} \in F_1$ be the sequence of terms
chosen in line~\ref{line:LiftReduce_choose} of \LiftReduce{} when this
algorithm is applied to a leading syzygy term $s \in \Lead_\succ(\Syz(G))$.
Then we have
\[
s \succ m_1 \gen_{i_1} \succ \ldots \succ m_k \gen_{i_k} \,.
\]
However, the terms $m \gen_i \in F_1$ chosen in \LiftHybrid{}, \LiftTree{}, and
\LiftSubtree{} satisfy $s \succ m \gen_i$, but they are not necessarily
ordered.
\end{rem}

\LiftTree{} has two main advantages in comparison to \LiftHybrid{}. First, no
reductions as in line~\ref{line:LiftHybrid_reduction} of \LiftHybrid{} occur.
Second, the results of \LiftSubtree{} can be cached and reused. We will see an
example for this in the next section.

\begin{rem}
It is worth mentioning that the proposed algorithm can be easily parallelized.
First of all, it is inherently parallel in several ways: Whenever a single
leading syzygy term is lifted to a syzygy by \LiftTree{}, the different
branches of the resulting tree, which correspond to the recursive calls of
\LiftSubtree{}, are independent of each other and can thus be treated in
parallel. (Note that this approach should be implemented in such a way that it
works well with the caching of partial results.) Likewise, in the computation
of a syzygy module via Algorithm~\ref{alg:SyzLift}, the leading syzygy terms
can be treated in parallel. In view of a whole resolution, we may start to lift
the leading syzygy terms in the lower syzygy modules while we are still
computing the Schreyer frame, see Section~\ref{ssec:schreyer_frame}, for the
higher ones. Note that, however, the time to compute the Schreyer frame is
almost negligible in most cases.

Over the rational numbers, one could also apply modular methods to our
algorithm, that is, one could do the computation modulo several primes in
parallel and then lift the results back to characteristic zero. Both approaches
are, however, subject to future research.
\end{rem}

\section{Example}\label{sec:example}

In this section, we give an example in order to illustrate the differences
between the three approaches which we presented in the previous section. The
example has been chosen in such a way that it shows the benefits, but also
possible drawbacks of the new methods. However, note that a considerable
speed-up can only be expected for large examples.

Throughout this section, let $F_0 := R := \Q[w,x,y,z]$ be endowed with the
lexicographical ordering, denoted by $>$. We compute the first syzygy module of
$G := (f_1, f_2, f_3) \subset R$ with
\begin{align*}
f_1 &:= wx+wz+x^2+2xz-z^2 \,, \\
f_2 &:= wy-wz-xz-yz-2z^2 \,, \\
f_3 &:= xy+z^2 \,.
\end{align*}
Note that $G$ is a Gr\"obner basis w.r.t.\@ $>$. Let $\succ$ be the Schreyer
ordering on $F_1 := R^3$ induced by $>$ and $G$. We can use
Algorithm~\ref{alg:LeadSyz} to check that a minimal generating set of the
leading syzygy module of $G$ w.r.t.\@ $\succ$ is given by
\[
\mcL = \{ x \cdot \gen_2,\, w \cdot \gen_3 \} \subset F_1 \,.
\]

\begin{figure}[tb]
\caption[\textsc{LiftReduce}/\textsc{LiftHybrid} applied to
$x \cdot \gen_2$ and $w \cdot \gen_3$]%
{\textsc{LiftReduce}/\textsc{LiftHybrid} applied to \\
$x \cdot \gen_2$ and $w \cdot \gen_3$}%
\label{fig:LiftReduce}
\begin{center}
\small
\hfill
\begin{tikzpicture}[syzdiagram]
\node { $\begin{aligned}
\bar{s}_1 &:= x \cdot \gen_2 \\
        g &:= wxy -wxz \lot{{}-x^2z} -xyz \lot{{}-2xz^2}
\end{aligned}$ }
child { node { $\begin{aligned}
\bar{s}_1 &:= \bar{s}_1 -y \cdot \gen_1 \\
        g &:= -wxz -wyz -x^2y \lot{{}-x^2z} \\
          &\phantom{{}:={}} -3xyz \lot{{}-2xz^2 +yz^2}
\end{aligned}$ }
child { node { $\begin{aligned}
\bar{s}_1 &:= \bar{s}_1 +z \cdot \gen_1 \\
        g &:= -wyz \lot{{}+wz^2} -x^2y -3xyz \lot{{}+yz^2 -z^3}
\end{aligned}$ }
child { node { $\begin{aligned}
\bar{s}_1 &:= \bar{s}_1 +z \cdot \gen_2 \\
        g &:= -x^2y -3xyz \lot{{}-xz^2 -3z^3}
\end{aligned}$ }
child { node { $\begin{aligned}
\bar{s}_1 &:= \bar{s}_1 +x \cdot \gen_3 \\
        g &:= -3xyz \lot{{}-3z^3}
\end{aligned}$ }
child { node { $\begin{aligned}
\bar{s}_1 &:= \bar{s}_1 +3z \cdot \gen_3 \\
        g &:= 0
\end{aligned}$ } } } } } };
\end{tikzpicture}
\hfill
\begin{tikzpicture}[syzdiagram]
\node { $\begin{aligned}
\bar{s}_2 &:= w \cdot \gen_3 \\
        g &:= wxy \lot{{}+wz^2}
\end{aligned}$ }
child { node { $\begin{aligned}
\bar{s}_2 &:= \bar{s}_2 -y \cdot \gen_1 \\
        g &:= -wyz \lot{{}+wz^2} -x^2y \\
          &\phantom{{}:={}} -2xyz \lot{{}+yz^2}
\end{aligned}$ }
child { node { $\begin{aligned}
\bar{s}_2 &:= \bar{s}_2 +z \cdot \gen_2 \\
        g &:= -x^2y -2xyz \lot{{}-xz^2 -2z^3}
\end{aligned}$ }
child { node { $\begin{aligned}
\bar{s}_2 &:= \bar{s}_2 +x \cdot \gen_3 \\
        g &:= -2xyz \lot{{}-2z^3}
\end{aligned}$ }
child { node { $\begin{aligned}
\bar{s}_2 &:= \bar{s}_2 +2z \cdot \gen_3 \\
        g &:= 0
\end{aligned}$ } } } } };
\end{tikzpicture}
\hfill
\end{center}
\end{figure}
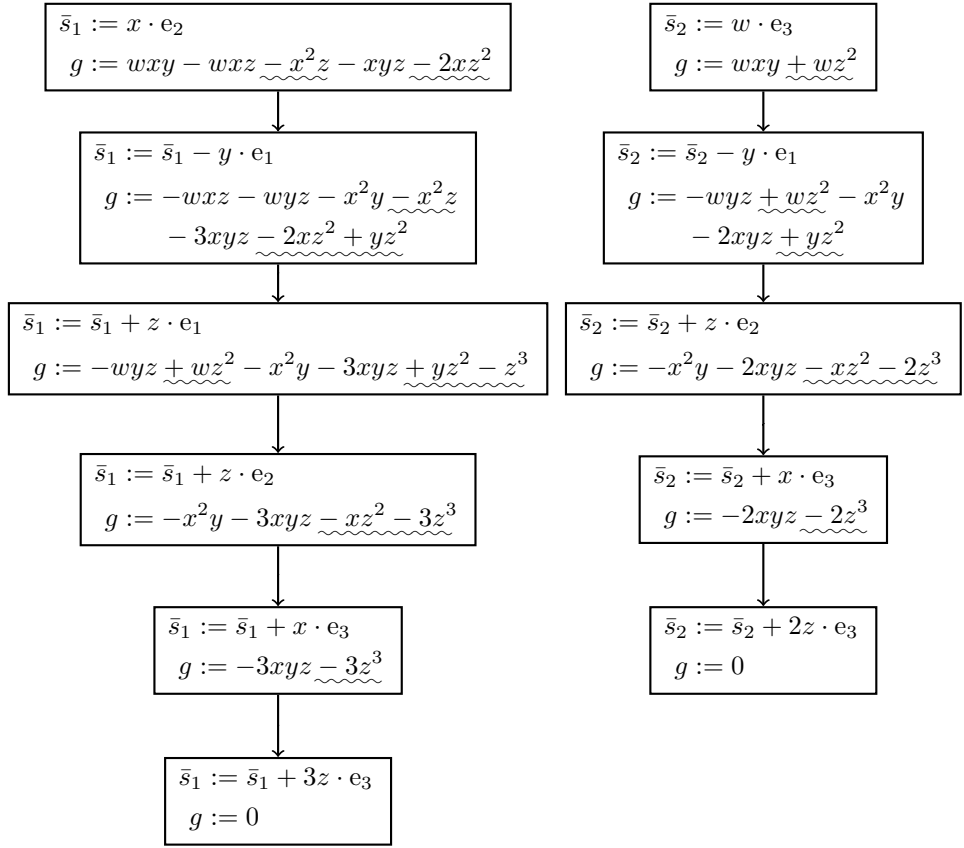

Our goal is to extend these leading syzygy terms to generators of the syzygy
module. Let us first consider the usual \textsc{LiftReduce} approach
(Algorithm~\ref{alg:LiftReduce}). Flow charts of \textsc{LiftReduce} applied to
the two leading syzygy terms above are shown in Figure~\ref{fig:LiftReduce}.

They start with the input term on the syzygy level and its image $g$ under
$\psi: F_1 \rightarrow F_0$, $\gen_i \mapsto f_i$, on the level of $F_0$. At
each step, the leading term of $g$ is reduced w.r.t.\@ $G$ while $\bar{s}_1$
and $\bar{s}_2$ keep track of these reductions. Both charts have the shape of a
chain because every step depends on the previous one. We could choose a
different reduction at the first step of the diagram on the right hand side,
but there is no other choice at the other steps. The process ends when $g = 0$
is reached, and we finally get the syzygies
\begin{align*}
\bar{s}_1 &= (-y+z) \cdot \gen_1 + (x+z) \cdot \gen_2 + (x+3z) \cdot \gen_3
\in \Syz(G) \; \text{ and} \\
\bar{s}_2 &= -y \cdot \gen_1 + z \cdot \gen_2 + (w+x+2z) \cdot \gen_3
\in \Syz(G)
\end{align*}
as liftings of the leading syzygy terms $x \cdot \gen_2$ and $w \cdot \gen_3$,
respectively.

The main innovation of the algorithm \textsc{LiftHybrid}
(Algorithm~\ref{alg:LiftHybrid}) is to leave out the lower order terms in $g$.
This in turn allows us to choose, at each step, any of the remaining terms for
reduction, in contrast to \LiftReduce{}. Hence the terms in $g$ do not have to
be ordered at all. If we always choose, however, to reduce the leading term as
in \textsc{LiftReduce}, then the flow charts of \textsc{LiftHybrid} applied to
$x \cdot \gen_2$ and $w \cdot \gen_3$, respectively, can be obtained from those
for \textsc{LiftReduce} by leaving out the underlined lower order terms, see
Figure~\ref{fig:LiftReduce}.

\begin{figure}[tb]
\caption{\textsc{LiftTree} applied to $w \cdot \gen_3$}
\label{fig:LiftTree2}
\begin{center}
\small
\begin{tikzpicture}[syzdiagram]
\node { $\begin{aligned}
\bar{s}_2 &:= w \cdot \gen_3 \\
        T &:= \{ wxy, \lot{wz^2} \}
\end{aligned}$ }
child { node {cached:\\[1ex] $\begin{aligned}
\bar{s}_2 &:= -y \cdot \gen_1 +z \cdot \gen_2 +x \cdot \gen_3
              +2z \cdot \gen_3 \\
        T &:= \varnothing
\end{aligned}$ } };
\end{tikzpicture}
\end{center}
\end{figure}
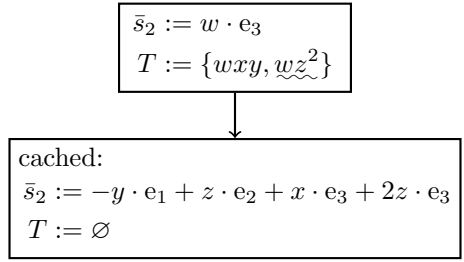

In \textsc{LiftTree} (Algorithm~\ref{alg:LiftTree}), the polynomial $g$ is
replaced by a set of terms denoted by $T$ and each term is treated
independently. The corresponding flow charts in Figure~\ref{fig:LiftTree1} and
Figure~\ref{fig:LiftTree2} thus have a tree structure where each node
represents one of the recursive calls of \textsc{LiftTree} and
\textsc{LiftSubtree}. In Figure~\ref{fig:LiftTree2}, the result of
$\LiftTree(x \cdot \gen_2)$ can be read off as the sum of all the terms
$\bar{s}_1$. Similarly, $\LiftTree(w \cdot \gen_3)$ yields
$-y \cdot \gen_1 + z \cdot \gen_2 + (w+x+2z) \cdot \gen_3$.

\begin{figure}[tb]
\caption{\textsc{LiftTree} applied to $x \cdot \gen_2$}
\label{fig:LiftTree1}
\begin{center}
\small
\begin{tikzpicture}[syzdiagram]
\node (1) { $\begin{aligned}
\bar{s}_1 &:= x \cdot \gen_2 \\
        T &:= \{ wxy, -wxz, \lot{-x^2z}, -xyz, \lot{-2xz^2} \}
\end{aligned}$ }
child { node[shift=(0:-2.1cm)] (1-1) { $\begin{aligned}
\bar{s}_1 &:= -y \cdot \gen_1 \\
        T &:= \{ -wyz, -x^2y, -2xyz, \lot{yz^2} \}
\end{aligned}$ }
child { node[shift=(0:1.5cm)] (1-1-1) { $\begin{aligned}
\bar{s}_1 &:= z \cdot \gen_2 \\
        T &:= \{ \lot{-wz^2, -xz^2, -yz^2, -2z^3} \}
\end{aligned}$ } }
child { node[rightof=1-1-1] (1-1-2) { $\begin{aligned}
\bar{s}_1 &:= x \cdot \gen_3 \\
        T &:= \{ \lot{xz^2} \}
\end{aligned}$ } }
child { node[rightof=1-1-2] (1-1-3) { $\begin{aligned}
\bar{s}_1 &:= 2z \cdot \gen_3 \\
        T &:= \{ \lot{2z^3} \}
\end{aligned}$ } } }
child { node[rightof=1-1] (1-2) { $\begin{aligned}
\bar{s}_1 &:= z \cdot \gen_1 \\
        T &:= \{ \lot{wz^2, x^2z, 2xz^2, -z^3} \}
\end{aligned}$ } }
child { node[rightof=1-2] (1-3) { $\begin{aligned}
\bar{s}_1 &:= z \cdot \gen_3 \\
        T &:= \{ \lot{z^3} \}
\end{aligned}$ } };
\end{tikzpicture}
\end{center}
\end{figure}
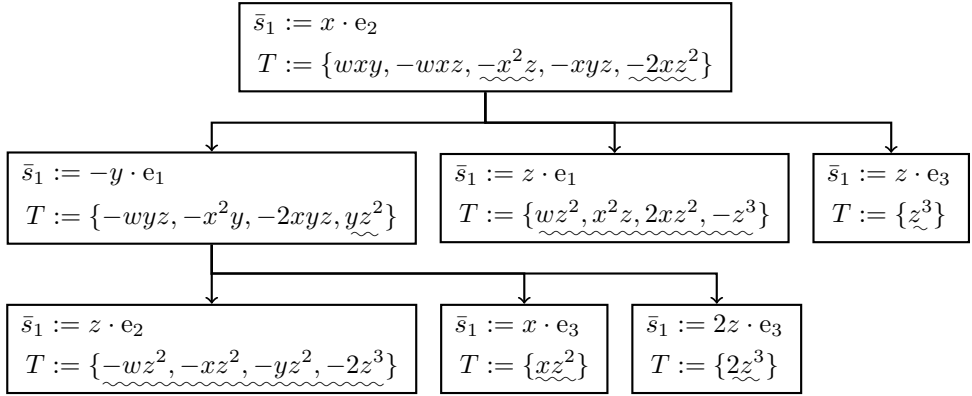

Again, the underlined lower order terms are left out. The process ends when
$T = \varnothing$ is reached in every branch. It is worth noting that although
each step resembles a reduction step, no reductions as in the first two
approaches occur. The main advantage of the \textsc{LiftTree} approach is that
intermediate results can be cached and reused. In Figure~\ref{fig:LiftTree2},
the term $wxy$ occurs as an element of $T$, but the whole subtree which
corresponds to this element has already been computed when \textsc{LiftTree}
was applied to $x \cdot \gen_2$ in Figure~\ref{fig:LiftTree1} and we can
therefore just plug in the cached result.

A possible drawback of this method can be observed in
Figure~\ref{fig:LiftTree1}: Two steps are necessary to compute the term
$(3z \cdot \gen_3)$ in the result whereas \textsc{LiftReduce} and
\textsc{LiftHybrid} need only one step for this. On the other hand, we could
also cache the result of $\LiftSubtree(z \cdot \gen_3)$ and reuse it for the
computation of $\LiftSubtree(2z \cdot \gen_3)$, of course.

As mentioned above, the leading terms of the two computed syzygies
\begin{align*}
\bar{s}_1 &= (-y+z) \cdot \gen_1 + (x+z) \cdot \gen_2 + (x+3z) \cdot \gen_3
\in \Syz(G) \; \text{ and} \\
\bar{s}_2 &= -y \cdot \gen_1 + z \cdot \gen_2 + (w+x+2z) \cdot \gen_3
\in \Syz(G)
\end{align*}
w.r.t.\@ the Schreyer ordering $\succ$ are $x \cdot \gen_2$ and
$w \cdot \gen_3$, respectively. Since these two terms belong to different
module components, the S-vector $S(\bar{s}_1, \bar{s}_2)$ of the two syzygies
is $0$. This implies $\Syz(\Syz(G)) = 0$ by Theorem~\ref{thm:Schreyer} applied
to $\Syz(G)$. We therefore get
\[
\mathcal{F}: \quad
0 \longrightarrow R^2
\overset{\phi_2}{\longrightarrow} R^3
\overset{\phi_1}{\longrightarrow} R
\longrightarrow R / \langle G \rangle_R
\longrightarrow 0
\]
as a free resolution of $R / \langle G \rangle_R$, where the maps are given by
$\phi_1(\gen_i) := f_i$, for $i = 1,2,3$, and $\phi_2(\gen_j) := \bar{s}_j$,
for $j = 1,2$. Since $G$ and $\{\bar{s}_1, \bar{s}_2 \}$ are reduced Gr\"obner
bases w.r.t.\@ the monomial orderings $>$ and $\succ$, respectively,
$\mathcal{F}$ is even a minimal free resolution.

\section{Timings and Statistics}\label{sec:timings}

In this section, we illustrate, by a number of examples, the speed-up achieved
by \LiftTree{} in comparison to other algorithms, and we also give detailed
statistics on these computations. In Subsection~\ref{ssec:agr}, we consider two
series of Artinian graded Gorenstein rings while Subsection~\ref{ssec:cnc_pcnc}
is devoted to randomly constructed canonical and Prym canonical nodal curves.

For each example, we compare the timings to compute a free resolution by
the Macaulay2 \cite{M2} command \verb+res+, by the \Singular{} \cite{DGPS}
command \verb+lres()+, and by the implementation of \LiftTree{} in the
\Singular{} library \texttt{schreyer.lib} \cite{M14}. Note that the first two
yield \emph{minimal} free resolution whereas the resolutions computed by
\LiftTree{} are in general \emph{non-minimal}. We therefore additionally
present the timings to compute the minimal Betti numbers via the \Singular{}
command \verb+betti()+ and to minimize the whole resolution via the \Singular{}
command \verb+minres()+, starting from the output of \LiftTree{} in both cases.
In separate tables, we also list, for each example, the number of terms in the
resolution computed by \LiftTree{} (excluding the first map which is given by
the input ideal) as well as the number of multiplications, additions, and
cancellations (that is, additions to zero) of coefficients in the ground field
which were performed during this computation. To measure the sparseness of the
resolution, the quotient of the number of terms divided by the number of
(matrix) entries $Q_\text{sparse}$ is also given. Finally, we present the
(minimal or non-minimal) Betti tables for selected examples.

In each case, we used the degree reverse lexicographic monomial ordering
(\texttt{dp} in \textsc{Singular}) and we computed a reduced Gr\"obner basis of
the input ideal w.r.t.\@ this ordering beforehand. The timings were computed on
an Intel Core i7-860 machine with 16 GB RAM and 4 physical (8 virtual) cores,
each with 2.8 GHz, running Fedora 20 (Linux kernel version 3.17.4). A dash
($-$) indicates that the computation did not finish within 24~hours.

\subsection{Artinian graded Gorenstein rings}%
\label{ssec:agr}

Our first series of examples consists of Artinian graded Gorenstein rings (AGR)
in $n+1$ variables of socle degree $d$ (see Tables~\ref{tab:timings_agr} to
\ref{tab:statistics_agr2}). We regard this family of zero-dimensional ideals as
a family where we can easily alter the number of variables, the degree, and the
sparseness. Of course it would be nice to produce higher dimensional random
examples, say, random ideals of codimension $c < n$ generated by $s$ forms of
degree $d_1, \ldots, d_s$ with $s \neq c$. However, since the Hilbert schemes
of such examples are not unirational in most cases, we do not know how to
construct such examples.

Artinian graded Gorenstein rings in $n+1$ variables of socle degree $d$ can be
obtained from a homogeneous form $f$ on $\PP^n$ of degree $d$ via apolarity,
see, for example, \cite{RS}. If we write $f = \ell_1^d + \ldots + \ell_s^d$ as
a sum of $s$ powers of linear forms, then the minimal such $s$ is called the
Waring rank of $f$. The shape of the minimal resolution and even the length of
the Gorenstein ring depends on $s$, and we invite the reader to prove in our
examples that for small $s$, the Waring decomposition is unique up to scalar, a
topic already studied for binary forms by Sylvester and recently picked up
again in \cite{OO}.

Below we report the findings for random examples over the finite field with
10,007 elements. It requires 24.12 seconds to compute the minimal resolution
for $s = 48$, while it just takes 3.12 seconds to get the minimal Betti
numbers.

\begin{figure}
\caption{Grayscale image of the $756 \times 840$ matrix from the non-minimal
resolution of the AGR example with $n = 6$, $d = 5$, and $s = 42$. Zero entries
are white, entries with one term gray, and entries with two terms black. There
are no entries with more than two terms.}%
\label{fig:agr_s42_m4}
\includegraphics[width=\textwidth]{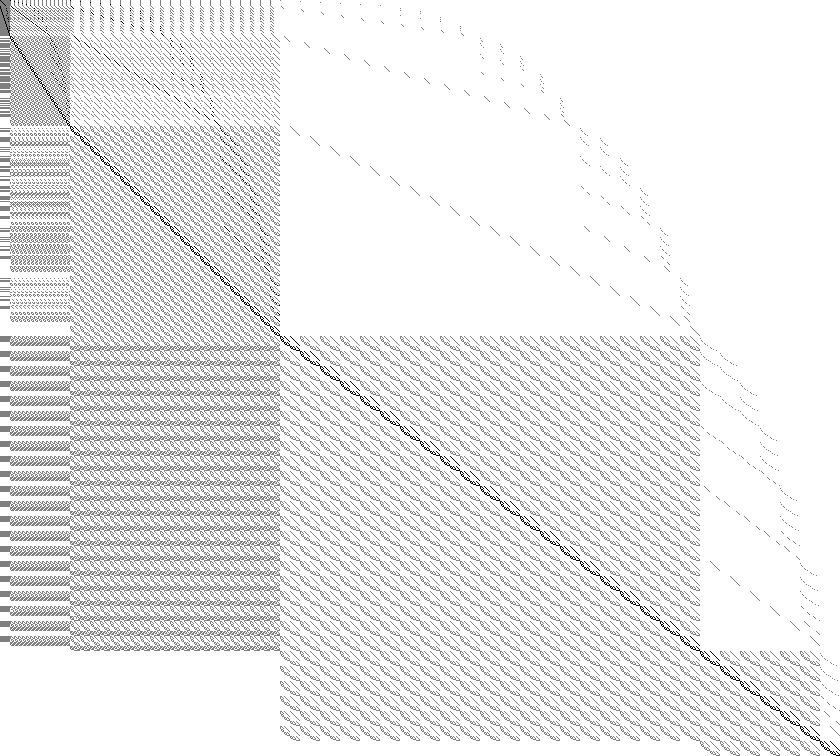}
\end{figure}

The total number of additions in most of these examples is smaller than the
total number of terms, and the number of cancellations (additions to zero) is
even much smaller. Furthermore, the matrices in the non-minimal resolution have
basically monomial entries as Figure~\ref{fig:agr_s42_m4} indicates. This
suggests that the algorithm \LiftTree{} to compute a non-minimal resolution is
numerically stable.

\begin{table}[htbp]
\centering
\caption{Timings for the AGR examples with $d = 5$ and $n = 6$ (in sec.)}%
\label{tab:timings_agr}

\begin{tabular}{|r||r|r|r|r|r|}
\hline
& \multicolumn{1}{c|}{Macaulay2} & \multicolumn{4}{c|}{\Singular{}} \\
\cline{2-6}
$s$ & \texttt{res} & \texttt{lres()} & \LiftTree & \texttt{betti()} &
    \texttt{minres()} \\
\hline\hline
$12$ &  5.28 &  6.35 &  0.06 & 0.02 &  0.96 \\ \hline
$18$ & 10.65 &  3.63 &  0.27 & 0.06 &  4.53 \\ \hline
$24$ & 26.89 &  7.59 &  1.06 & 0.19 & 33.31 \\ \hline
$30$ & 16.77 &  5.06 &  0.94 & 1.61 & 39.26 \\ \hline
$36$ & 44.45 & 28.16 &  0.73 & 2.09 & 11.74 \\ \hline
$42$ & 87.65 & 48.85 &  0.73 & 2.38 & 23.39 \\ \hline
$48$ & 87.61 & 48.69 &  0.73 & 2.39 & 23.39 \\ \hline
\end{tabular}
\end{table}

\begin{table}
\centering
\caption{Statistics for the non-minimal resolutions of the AGR examples with
$d = 5$ and $n = 6$}%
\label{tab:statistics_agr}

\begin{tabular}{|r||r||r|r|r|r|r|}
\hline
$s$ & sec. & \#Terms & \#Mult. & \#Add. & \#Canc. & $Q_{\text{sparse}}$ \\
\hline\hline
$12$ & 0.06 &  34,963 &    74,190 &    38,312 &  2,163 & 0.155 \\ \hline
$18$ & 0.27 & 123,144 &   700,889 &   573,143 &  7,267 & 0.224 \\ \hline
$24$ & 1.06 & 316,492 & 5,961,627 & 5,638,864 & 16,813 & 0.276 \\ \hline
$30$ & 0.94 & 319,580 &   627,508 &   315,310 &  2,496 & 0.177 \\ \hline
$36$ & 0.73 & 294,730 &   447,245 &   162,996 &    331 & 0.195 \\ \hline
$42$ & 0.73 & 294,762 &   447,249 &   163,002 &    324 & 0.195 \\ \hline
$48$ & 0.73 & 294,746 &   447,260 &   162,992 &    334 & 0.195 \\ \hline
\end{tabular}
\end{table}

\begin{table}
\centering
\caption{Minimal Betti table for the AGR example with $n = 6$, $d = 5$, and
$s = 18$}%
\label{tab:betti_agr_n6d5s18}

\begin{tabular}{rrrrrrrrr}
       &    0 &    1 &    2 &    3 &    4 &    5 &    6 &    7 \\
\hline
    0: &    1 &    - &    - &    - &    - &    - &    - &    - \\
    1: &    - &   10 &    4 &    - &    - &    - &    - &    - \\
    2: &    - &    - &   60 &  136 &  130 &   60 &   11 &    - \\
    3: &    - &   11 &   60 &  130 &  136 &   60 &    - &    - \\
    4: &    - &    - &    - &    - &    - &    4 &   10 &    - \\
    5: &    - &    - &    - &    - &    - &    - &    - &    1 \\
\hline
total: &    1 &   21 &  124 &  266 &  266 &  124 &   21 &    1 \\
\end{tabular}
\end{table}

\begin{table}
\centering
\caption{Minimal Betti table for the AGR example with $n = 6$, $d = 5$, and
$s = 24$}%
\label{tab:betti_agr_n6d5s24}

\begin{tabular}{rrrrrrrrr}
       &    0 &    1 &    2 &    3 &    4 &    5 &    6 &    7 \\
\hline
    0: &    1 &    - &    - &    - &    - &    - &    - &    - \\
    1: &    - &    4 &    - &    - &    - &    - &    - &    - \\
    2: &    - &   32 &  150 &  256 &  220 &   96 &   17 &    - \\
    3: &    - &   17 &   96 &  220 &  256 &  150 &   32 &    - \\
    4: &    - &    - &    - &    - &    - &    - &    4 &    - \\
    5: &    - &    - &    - &    - &    - &    - &    - &    1 \\
\hline
total: &    1 &   53 &  246 &  476 &  476 &  246 &   53 &    1 \\
\end{tabular}
\end{table}

\begin{table}
\centering
\caption{Minimal Betti table for the AGR examples with $n = 6$, $d = 5$, and
$s \geq 42$}%
\label{tab:betti_agr_n6d5s42}

\begin{tabular}{rrrrrrrrr}
       &    0 &    1 &    2 &    3 &    4 &    5 &    6 &    7 \\
\hline
    0: &    1 &    - &    - &    - &    - &    - &    - &    - \\
    1: &    - &    - &    - &    - &    - &    - &    - &    - \\
    2: &    - &   56 &  189 &  216 &    - &    - &    - &    - \\
    3: &    - &    - &    - &    - &  216 &  189 &   56 &    - \\
    4: &    - &    - &    - &    - &    - &    - &    - &    - \\
    5: &    - &    - &    - &    - &    - &    - &    - &    1 \\
\hline
total: &    1 &   56 &  189 &  216 &  216 &  189 &   56 &    1 \\
\end{tabular}
\end{table}

\begin{table}
\centering
\caption{The minimal Betti tables for the AGR examples with $n = 6$, $d = 5$,
and $s \geq 28$ are the sum of Table~\ref{tab:betti_agr_n6d5s42} and the table
below.}%
\label{tab:betti_agr_n6d5s28}

\begin{tabular}{ccccccccc}
       &    0 &    1 &    2 &    3 &    4 &    5 &    6 &    7 \\
\hline
    0: &    - &    - &    - &    - &    - &    - &    - &    - \\
    1: &    - &    - &    - &    - &    - &    - &    - &    - \\
    2: &    - &    - & $\alpha(s)$ & $\beta(s)$ & $\gamma(s)$ &
                                  $\beta(s)$ & $\alpha(s)$ & - \\
    3: &    - & $\alpha(s)$ & $\beta(s)$ & $\gamma(s)$ &
                             $\beta(s)$ & $\alpha(s)$ & -  & - \\
    4: &    - &    - &    - &    - &    - &    - &    - &    - \\
    5: &    - &    - &    - &    - &    - &    - &    - &    - \\
\hline
\end{tabular} \\
\vspace{1ex}
where
$\alpha(s) := \max \{ 0,\, 6(31.5-s) \}$,\,
$\beta(s)  := \max \{ 0,\, 15(36-s) \}$,\,
$\gamma(s) := \max \{ 0,\, 20(42-s) \}$.
\end{table}

\begin{table}
\centering
\caption{Non-minimal Betti table for the AGR examples with $n = 6$, $d = 5$,
and $s \geq 34$}%
\label{tab:betti_agr_n6d5s42_nonmin}

\begin{tabular}{rrrrrrrrr}
       &    0 &    1 &    2 &    3 &    4 &    5 &    6 &    7 \\
\hline
    0: &    1 &    - &    - &    - &    - &    - &    - &    - \\
    1: &    - &    - &    - &    - &    - &    - &    - &    - \\
    2: &    - &   56 &  210 &  336 &  280 &  120 &   21 &    - \\
    3: &    - &   21 &  126 &  315 &  420 &  315 &  126 &   21 \\
    4: &    - &    6 &   36 &   90 &  120 &   90 &   36 &    6 \\
    5: &    - &    1 &    6 &   15 &   20 &   15 &    6 &    1 \\
\hline
total: &    1 &   84 &  378 &  756 &  840 &  540 &  189 &   28 \\
\end{tabular}
\end{table}

\begin{table}
\centering
\caption{Timings for the AGR examples with $d = 5$ and large~$s$ (in sec.)}%
\label{tab:timings_agr2}

\begin{tabular}{|r||r|r|r|r|r|}
\hline
& \multicolumn{1}{c|}{Macaulay2} & \multicolumn{4}{c|}{\Singular{}} \\
\cline{2-6}
$n$ & \texttt{res} & \texttt{lres()} & \LiftTree & \texttt{betti()} &
    \texttt{minres()} \\
\hline\hline
 $5$ &      2.47 &      1.17 &      0.06 &      0.06 &     0.49 \\ \hline
 $6$ &     87.74 &     48.99 &      0.73 &      2.38 &    23.39 \\ \hline
 $7$ &  1,782.44 &  1,027.15 &      9.97 &     60.14 &   322.45 \\ \hline
 $8$ & 50,203.20 & 32,647.20 &    137.25 &  1,114.14 & 9,002.31 \\ \hline
 $9$ &       $-$ &       $-$ &  1,852.42 & 16,100.67 &      $-$ \\ \hline
$10$ &       $-$ &       $-$ & 24,281.08 &       $-$ &      $-$ \\ \hline
\end{tabular}
\end{table}

\begin{table}
\centering
\caption{Statistics for the non-minimal resolutions of the AGR examples with
$d = 5$ and large $s$}%
\label{tab:statistics_agr2}

\begin{tabular}{|r||r||r|r|r|r|r|}
\hline
$n$ & sec. & \#Terms & \#Mult. & \#Add. & \#Canc. & $Q_{\text{sparse}}$ \\
\hline\hline
 $5$ &      0.06 &     59,903 &    101,264 &    44,790 &     63 & 0.282 \\
\hline
 $6$ &      0.73 &    294,746 &    447,232 &   162,989 &    327 & 0.195 \\
\hline
 $7$ &      9.97 &  1,292,567 &  1,761,229 &   496,922 &  1,185 & 0.130 \\
\hline
 $8$ &    137.25 &  5,179,579 &  6,433,983 & 1,323,234 &  3,562 & 0.084 \\
\hline
 $9$ &  1,852.42 & 19,311,659 & 22,322,538 & 3,166,754 &  9,013 & 0.053 \\
\hline
$10$ & 24,281.08 & 67,915,012 & 74,543,926 & 6,963,586 & 20,518 & 0.033 \\
\hline
\end{tabular}
\end{table}

\clearpage

\subsection{Canonical Nodal Curves and Prym Canonical Nodal Curves}%
\label{ssec:cnc_pcnc}

Our next two series of examples are randomly constructed canonical nodal curves
(CNC) (see Tables~\ref{tab:timings_cnc} and \ref{tab:statistics_cnc}) and Prym
canonical nodal curves (PCNC) (see Tables~\ref{tab:timings_pcnc} and
\ref{tab:statistics_pcnc}) of genus $g$. For the series of CNC examples, we
work over the finite field with 32,003 elements, whereas we consider the PCNC
examples over fields of random positive characteristic between 10,000 and
30,000 as listed in the second column of Table~\ref{tab:timings_pcnc}.

Canonical curves are a widely studied topic, see, for example, \cite{BS} for
some references. Our interest in the minimal resolution of Prym canonical
curves comes from \cite{CEFS}, where an analogue of Green's conjecture for Prym
curves was formulated. In terms of the resolutions computed with our method,
both conjectures say that the degree zero parts of the syzygy matrices in the
resolution have maximal rank. Computational results presented in \cite{CEFS}
indicate that the Prym-Green conjecture is very likely false for genus $g = 8$
and $g = 16$. It would be interesting to check experimentally whether $g = 24$
is an exception as well. But our timings suggest that this is way out of reach
even for our improved method. In fact, computing a resolution with \LiftTree{}
for the case $g = 24$ would take, as a very rough estimate, about 100,000~years
and 100~TB of memory.

\begin{figure}
\caption{Black-and-white image of the $1932 \times 1932$ constant submatrix
from the non-minimal resolution of the PCNC example with $g = 14$. Zero entries
are white, constant entries are black. Note the somehow fractal structure.}%
\label{fig:pcnc_g14_m5_submat}
\includegraphics[width=\textwidth]{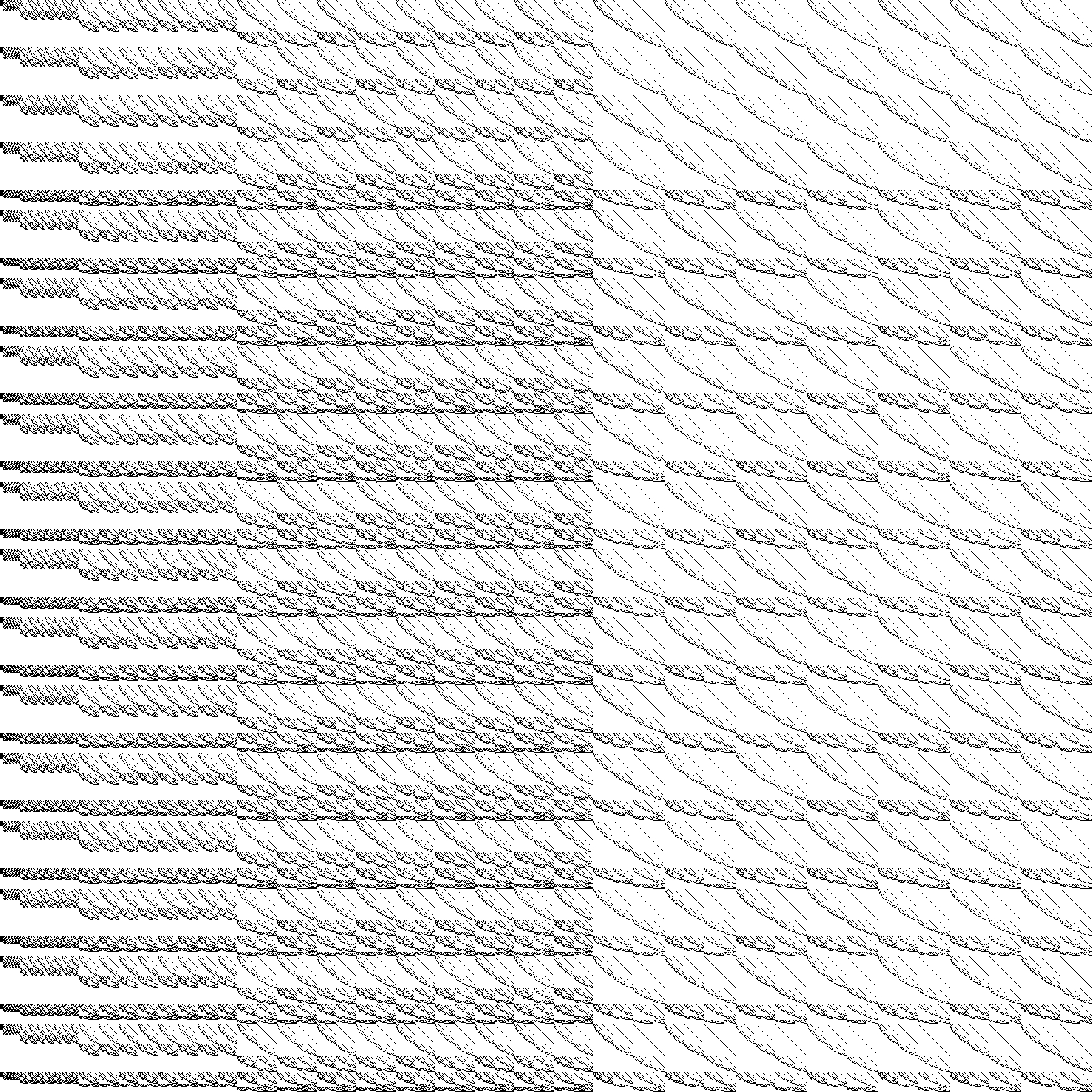}
\end{figure}

\begin{rem}
The non-minimal Betti numbers for the PCNC example with $g = 14$ as computed
with \LiftTree{} are shown in Table~\ref{tab:betti_pcnc_g14_nonmin}. Note that
the fifth differential $\phi_5$ contains a $1932 \times 1932$ submatrix with
constant entries. A black-and-white image of this submatrix is shown in
Figure~\ref{fig:pcnc_g14_m5_submat}. Using specialized software such as
FFLAS/FFPACK \cite{DGP} which is easily available via Sage \cite{Sage},
checking that this matrix has full rank takes only 0.92~seconds. Thus,
including the time to compute the non-minimal resolution (224.15~seconds, cf.\@
Table~\ref{tab:timings_pcnc}), the verification of the Prym-Green conjecture
for this case takes 225.07~seconds using our new method. If we only consider
the time to compute the non-minimal resolution up to the fifth differential
(129.87~seconds, cf.\@ Table~\ref{tab:timings_pcnc_g14_nonmin}), this can be
reduced to 130.79~seconds.
\end{rem}

\begin{rem}\label{rem:pcnc_g16}
For the PCNC example with $g = 16$, the non-minimal Betti numbers as computed
with \LiftTree{} are shown in Table~\ref{tab:betti_pcnc_g16_nonmin}. Here, the
sixth differential $\phi_6$ contains a $8910 \times 8910$ submatrix with
constant entries. Using FFLAS/FFPACK \cite{DGP}, it takes 55.24~seconds to
check that the kernel of this submatrix is one-dimensional. Including the time
to compute the non-minimal resolution up to the sixth differential
(10,776.54~seconds, cf.\@ Table~\ref{tab:timings_pcnc_g16_nonmin}), it thus
takes little more than 3~hours to check that the verification of the Prym-Green
conjecture based on nodal curves fails in this case. This is a substantial
improvement of the running time for this problem in comparison to the time
needed in \cite{CEFS}.

We also managed to compute the minimal Betti numbers for this example using the
\textsc{Singular} command \verb+betti()+, cf.\@ Table~\ref{tab:betti_pcnc_g16}.
This computation took more than 37~hours, but it was run on a different compute
server than the other examples due to the memory consumption exceeding 16~GB.

Finally, starting from the non-minimal resolution, we even succeeded to compute
the syzygy scheme of the Prym-Green extra syzygy which prevents the minimal
resolution of the curve from being pure, cf.\@ Table~\ref{tab:betti_pcnc_g16}.
However, we found that in our example, this syzygy scheme is identical to the
curve itself and does therefore not contribute any geometric information.
\end{rem}

\begin{table}[htbp]
\centering
\caption{Timings for the CNC examples (in sec.)}%
\label{tab:timings_cnc}

\begin{tabular}{|r||r|r|r|r|r|}
\hline
& \multicolumn{1}{c|}{Macaulay2} & \multicolumn{4}{c|}{\Singular{}} \\
\cline{2-6}
$g$ & \texttt{res} & \texttt{lres()} & \LiftTree & \texttt{betti()} &
    \texttt{minres()} \\
\hline\hline
 $5$ &        0 &         0 &         0 &        0 &         0 \\ \hline
 $6$ &     0.01 &         0 &         0 &        0 &         0 \\ \hline
 $7$ &     0.07 &      0.02 &      0.01 &        0 &         0 \\ \hline
 $8$ &     1.41 &      0.33 &      0.02 &        0 &      0.04 \\ \hline
 $9$ &    25.45 &      6.65 &      0.07 &     0.02 &      0.97 \\ \hline
$10$ &   452.10 &    113.48 &      0.31 &     0.38 &      9.28 \\ \hline
$11$ & 7,945.04 &  2,017.32 &      1.83 &     6.87 &    155.52 \\ \hline
$12$ &      $-$ & 30,495.55 &     15.58 &    89.66 &  1,238.28 \\ \hline
$13$ &      $-$ &       $-$ &    142.20 & 1,005.98 & 23,877.02 \\ \hline
$14$ &      $-$ &       $-$ &  1,351.59 & 9,645.67 &       $-$ \\ \hline
$15$ &      $-$ &       $-$ & 12,935.45 &      $-$ &       $-$ \\ \hline
\end{tabular}
\end{table}

\begin{table}[htbp]
\centering
\caption{Statistics for the non-minimal resolutions of the CNC examples}%
\label{tab:statistics_cnc}

\begin{tabular}{|r||r||r|r|r|r|r|}
\hline
$g$ & sec. & \#Terms & \#Mult. & \#Add. & \#Canc. & $Q_{\text{sparse}}$ \\
\hline\hline
 $5$ &         0 &        258 &        466 &     226 &     0 & 3.583 \\
\hline
 $6$ &         0 &      1,411 &      2,161 &     831 &     0 & 2.520 \\
\hline
 $7$ &      0.01 &      6,038 &      8,104 &   2,323 &     3 & 1.677 \\
\hline
 $8$ &      0.02 &     22,343 &     27,123 &   5,453 &    18 & 1.075 \\
\hline
 $9$ &      0.07 &     75,054 &     84,804 &  11,319 &    63 & 0.669 \\
\hline
$10$ &      0.31 &    235,179 &    253,212 &  21,434 &   169 & 0.408 \\
\hline
$11$ &      1.83 &    699,758 &    730,490 &  37,790 &   380 & 0.244 \\
\hline
$12$ &     15.58 &  1,998,583 &  2,047,201 &  62,944 &   758 & 0.144 \\
\hline
$13$ &    142.20 &  5,522,774 &  5,593,998 & 100,053 & 1,393 & 0.084 \\
\hline
$14$ &  1,351.59 & 14,854,811 & 14,949,653 & 152,971 & 2,382 & 0.049 \\
\hline
$15$ & 12,935.45 & 39,056,118 & 39,164,376 & 226,312 & 3,866 & 0.028 \\
\hline
\end{tabular}
\end{table}

\begin{table}
\centering
\caption{Timings for the PCNC examples (in sec.)}%
\label{tab:timings_pcnc}

\begin{tabular}{|r|r||r|r|r|r|r|}
\hline
&& \multicolumn{1}{c|}{Macaulay2} & \multicolumn{4}{c|}{\Singular{}} \\
\cline{3-7}
$g$ & Char. & \texttt{res} & \texttt{lres()} & \LiftTree & \texttt{betti()} &
    \texttt{minres()} \\
\hline\hline
 $6$ & 22,669 &         0 &        0 &         0 &         0 &         0 \\
\hline
 $7$ & 10,151 &      0.02 &     0.01 &         0 &         0 &         0 \\
\hline
 $8$ & 15,187 &      0.18 &     0.05 &      0.01 &         0 &      0.01 \\
\hline
 $9$ & 18,947 &      2.96 &     0.75 &      0.04 &         0 &      0.11 \\
\hline
$10$ & 13,523 &     64.78 &    19.56 &      0.15 &      0.02 &      1.96 \\
\hline
$11$ & 25,219 &    901.14 &   344.62 &      0.64 &      0.40 &     23.11 \\
\hline
$12$ & 11,777 & 16,597.90 & 6,261.13 &      3.43 &      6.97 &    305.45 \\
\hline
$13$ & 24,379 &       $-$ &      $-$ &     25.95 &     98.34 &  3,032.23 \\
\hline
$14$ & 16,183 &       $-$ &      $-$ &    224.15 &  1,084.88 & 40,106.95 \\
\hline
$15$ & 20,873 &       $-$ &      $-$ &  2,002.51 & 10,953.23 &       $-$ \\
\hline
$16$ & 12,451 &       $-$ &      $-$ & 18,612.82 &       $-$ &       $-$ \\
\hline
\end{tabular}
\end{table}

\begin{table}
\centering
\caption{Statistics for the non-minimal resolutions of the PCNC examples}%
\label{tab:statistics_pcnc}

\begin{tabular}{|r||r||r|r|r|r|r|}
\hline
$g$ & sec. & \#Terms & \#Mult. & \#Add. & \#Canc. & $Q_{\text{sparse}}$ \\
\hline\hline
 $6$ &         0 &        501 &        534 &        90 &     0 & 2.088 \\
\hline
 $7$ &         0 &      2,818 &      5,297 &     2,663 &     1 & 2.271 \\
\hline
 $8$ &      0.01 &     12,974 &     24,165 &    11,642 &     8 & 1.917 \\
\hline
 $9$ &      0.04 &     48,711 &     79,096 &    31,409 &    14 & 1.390 \\
\hline
$10$ &      0.15 &    165,346 &    232,455 &    69,389 &    43 & 0.945 \\
\hline
$11$ &      0.64 &    524,473 &    654,596 &   135,098 &   110 & 0.617 \\
\hline
$12$ &      3.43 &  1,582,334 &  1,812,443 &   240,603 &   285 & 0.392 \\
\hline
$13$ &     25.95 &  4,594,249 &  4,974,596 &   401,889 &   594 & 0.244 \\
\hline
$14$ &    224.15 & 12,931,450 & 13,525,642 &   637,340 & 1,154 & 0.149 \\
\hline
$15$ &  2,002.51 & 35,482,705 & 36,367,579 &   970,090 & 2,018 & 0.090 \\
\hline
$16$ & 18,612.82 & 95,281,070 & 96,541,345 & 1,427,327 & 3,409 & 0.054 \\
\hline
\end{tabular}
\end{table}

\begin{table}
\centering
\caption{Non-minimal Betti table for the PCNC example with $g = 14$}%
\label{tab:betti_pcnc_g14_nonmin}

\addtolength{\tabcolsep}{-1pt}
\begin{tabular}{rrrrrrrrrrrrr}
       & 0 &  1 &   2 &    3 &    4 &    5 &    6 &    7 &    8 &   9 &  10 &
  11 \\ \hline
0:     & 1 &  - &   - &    - &    - &    - &    - &    - &    - &   - &   - &
   - \\
1:     & - & 52 & 303 &  882 & 1596 & 1932 & 1602 &  903 &  332 &  72 &   7 &
   - \\
2:     & - & 17 & 167 &  738 & 1932 & 3318 & 3906 & 3192 & 1788 & 657 & 143 &
  14 \\ \hline
total: & 1 & 69 & 470 & 1620 & 3528 & 5250 & 5508 & 4095 & 2120 & 729 & 150 &
  14 \\
\end{tabular}
\end{table}

\begin{table}
\centering
\caption{Timings (in sec.) and statistics for the individual differentials
$\phi_i$ in the non-minimal resolution of the PCNC example with $g = 14$. Note
that $\phi_1$ is given by the input ideal.}%
\label{tab:timings_pcnc_g14_nonmin}

\begin{tabular}{|r||r|r|r|r|r|}
\hline
$i$ & \#Generators & \#Terms & $Q_\text{sparse}$ & Time & $\sim$ \#Terms/sec.
\\ \hline\hline
      1 &     69 &      3,064 & 44.406 &    $-$ &       $-$ \\ \hline
      2 &    470 &     50,208 &  1.548 &   0.54 &    92,978 \\ \hline
      3 &  1,620 &    371,814 &  0.488 &   6.25 &    59,490 \\ \hline
      4 &  3,528 &  1,290,516 &  0.226 &  37.98 &    33,979 \\ \hline
      5 &  5,250 &  2,639,132 &  0.142 &  85.10 &    31,012 \\ \hline
      6 &  5,508 &  3,436,908 &  0.119 &  70.07 &    49,050 \\ \hline
      7 &  4,095 &  2,917,668 &  0.129 &  21.21 &   137,561 \\ \hline
      8 &  2,120 &  1,593,830 &  0.184 &   2.72 &   585,967 \\ \hline
      9 &    729 &    530,548 &  0.343 &   0.25 & 2,122,192 \\ \hline
     10 &    150 &     94,488 &  0.864 &   0.02 & 4,724,400 \\ \hline
     11 &     14 &      6,338 &  3.018 &   0.01 & (633,800) \\ \hline
\hline
2 to 11 & 23,484 & 12,931,450 &  0.149 & 224.15 &    57.691 \\
\hline
\end{tabular}
\end{table}

\begin{table}
\centering
\caption{Non-minimal Betti table for the PCNC example with $g = 16$}%
\label{tab:betti_pcnc_g16_nonmin}

\begin{tabular}{rrrrrrrrrrr}
       & 0 &  1 &   2 &    3 &    4 &     5 &     6 &     7 &     8 & \\
\cline{1-10}\cdashline{11-11}
0:     & 1 &  - &   - &    - &    - &     - &     - &     - &     - & \\
1:     & - & 75 & 539 & 1980 & 4653 &  7590 &  8910 &  7623 &  4730 & \\
2:     & - & 19 & 225 & 1221 & 4015 &  8910 & 14058 & 16170 & 13662 & \\
\cline{1-10}\cdashline{11-11}
total: & 1 & 94 & 764 & 3201 & 8668 & 16500 & 22968 & 23793 & 18392 & \\
\end{tabular}

\vspace{2ex}
\begin{tabular}{rrrrrr}
\phantom{X} &     9 &   10 &   11 &  12 & 13 \\
\cdashline{1-1}\cline{2-6}
            &     - &    - &    - &   - &  - \\
            &  2079 &  615 &  110 &   9 &  - \\
            &  8415 & 3685 & 1089 & 195 & 16 \\
\cdashline{1-1}\cline{2-6}
            & 10494 & 4300 & 1199 & 204 & 16 \\
\end{tabular}
\end{table}

\begin{table}
\centering
\caption{Minimal Betti table for the PCNC example with $g = 16$}%
\label{tab:betti_pcnc_g16}

\begin{tabular}{rrrrrrrrrrr}
       & 0 &  1 &   2 &    3 &    4 &    5 &    6 &     7 &     8 & \\
\cline{1-10}\cdashline{11-11}
0:     & 1 &  - &   - &    - &    - &    - &    - &     - &     - & \\
1:     & - & 75 & 520 & 1755 & 3432 & 3575 &    1 &     - &     - & \\
2:     & - &  - &   - &    - &    - &    1 & 6435 & 11440 & 11583 & \\
\cline{1-10}\cdashline{11-11}
total: & 1 & 75 & 520 & 1755 & 3432 & 3576 & 6436 & 11440 & 11583 & \\
\end{tabular}

\vspace{2ex}
\begin{tabular}{rrrrrr}
\phantom{X} &    9 &   10 &   11 &  12 & 13 \\
\cdashline{1-1}\cline{2-6}
            &    - &    - &    - &   - &  - \\
            &    - &    - &    - &   - &  - \\
            & 7800 & 3575 & 1080 & 195 & 16 \\
\cdashline{1-1}\cline{2-6}
            & 7800 & 3575 & 1080 & 195 & 16 \\
\end{tabular}
\end{table}

\begin{table}
\centering
\caption{Timings (in sec.) and statistics for the individual differentials
$\phi_i$ in the non-minimal resolution of the PCNC example with $g = 16$. Note
that $\phi_1$ is given by the input ideal.}%
\label{tab:timings_pcnc_g16_nonmin}

\begin{tabular}{|r||r|r|r|r|r|}
\hline
$i$ & \#Generators & \#Terms & $Q_\text{sparse}$ & Time & $\sim$ \#Terms/sec.
\\
\hline\hline
      1 &      94 &      4,706 & 50.064 &       $-$ &       $-$ \\ \hline
      2 &     764 &     96,185 &  1.339 &      1.75 &    54,963 \\ \hline
      3 &   3,201 &    881,460 &  0.360 &     51.64 &    17,069 \\ \hline
      4 &   8,668 &  3,884,192 &  0.140 &    669.57 &     5,801 \\ \hline
      5 &  16,500 & 10,428,510 &  0.073 &  3,332.74 &     3,129 \\ \hline
      6 &  22,968 & 18,632,465 &  0.049 &  6,720.84 &     2,772 \\ \hline
      7 &  23,793 & 23,035,835 &  0.042 &  5,663.46 &     4,067 \\ \hline
      8 &  18,392 & 19,959,667 &  0.046 &  1,897.21 &    10,521 \\ \hline
      9 &  10,494 & 12,027,377 &  0.062 &    260.48 &    46,174 \\ \hline
     10 &   4,300 &  4,889,758 &  0.108 &     14.46 &   338,158 \\ \hline
     11 &   1,199 &  1,257,155 &  0.244 &      0.59 & 2,130,771 \\ \hline
     12 &     204 &    178,565 &  0.730 &      0.08 & 2,232,063 \\ \hline
     13 &      16 &      9,901 &  3.033 &         0 &       $-$ \\ \hline
\hline
2 to 13 & 110,499 & 95,281,070 &  0.054 & 18,612.82 &     5,119 \\
\hline
\end{tabular}
\end{table}

\clearpage

\section{Acknowledgements}

We would like to thank Wolfram Decker for his steady encouragement during this
project.

\end{document}